\documentclass[preprint,12pt]{amsart}

\usepackage[margin=1.35in]{geometry}

\makeatletter
\def\ps@pprintTitle{%
  \let\@oddhead\@empty
  \let\@evenhead\@empty
  \let\@oddfoot\@empty
  \let\@evenfoot\@oddfoot
}
\makeatother
\usepackage{amsmath,amssymb,amsthm,mathtools}
\usepackage{enumitem}
\usepackage{booktabs}
\usepackage{longtable}
\usepackage{array}
\usepackage[expansion=false]{microtype}
\usepackage[colorlinks=true,linkcolor=blue!60!black,citecolor=green!40!black,urlcolor=blue!60!black]{hyperref}
\usepackage{xcolor}

\numberwithin{equation}{section}
\theoremstyle{plain}
\newtheorem{theorem}{Theorem}[section]
\newtheorem{proposition}[theorem]{Proposition}
\newtheorem{lemma}[theorem]{Lemma}
\newtheorem{corollary}[theorem]{Corollary}
\theoremstyle{definition}

\newtheorem{remark}[theorem]{Remark}

\newcommand{\Z}{\mathbb{Z}}
\newcommand{\Q}{\mathbb{Q}}
\newcommand{\C}{\mathbb{C}}
\newcommand{\F}{\mathbb{F}}
\newcommand{\Qbar}{\overline{\mathbb{Q}}}

\newcommand{\h}{\mathbb{H}}
\newcommand{\PP}{\mathbb{P}}
\DeclareMathOperator{\SL}{SL}
\DeclareMathOperator{\GL}{GL}

\DeclareMathOperator{\End}{End}

\DeclareMathOperator{\Br}{Br}
\DeclareMathOperator{\Tr}{Tr}
\DeclareMathOperator{\cond}{cond}
\DeclareMathOperator{\lcm}{lcm}

\newcommand{\new}{\mathrm{new}}
\newcommand{\diam}[1]{\langle #1\rangle}
\newcommand{\mat}[4]{\begin{pmatrix}#1&#2\\#3&#4\end{pmatrix}}
\newcommand{\smat}[4]{\left(\begin{smallmatrix}#1&#2\\#3&#4\end{smallmatrix}\right)}

\begin{document}
\title{A completely decomposable modular Jacobian of dimension $1777$}

\author{Robert Auffarth}
\address{Departamento de Matem\'aticas, Facultad de Ciencias, Universidad de Chile, Santiago, Chile}
\email{rfauffar@uchile.cl}

\date{October 2026}

\begin{abstract}
We exhibit a congruence subgroup $\Gamma\subset\SL(2,\Z)$ of level $120$ and index $23040$ whose modular curve $X_\Gamma$ has genus $1777$ and whose Jacobian is isogenous over $\Qbar$ to a product of elliptic curves. This surpasses the previous record of a completely decomposable complex Jacobian of dimension $1297$ found by Ekedahl and Serre in $1993$.
\end{abstract}

\maketitle

\section{Introduction}

An abelian variety over an algebraically closed field is said to be \emph{completely decomposable} if it is isogenous to a product of elliptic curves. In 1993, Ekedahl and Serre \cite{ES93} posed two questions which are still open:
\begin{enumerate}[label=(Q\arabic*)]
\item Is there, for every $g\geq1$, a curve of genus $g$ with completely decomposable Jacobian?
\item Are the genera of such curves bounded?
\end{enumerate}

In their article, they constructed examples in many genera, both from modular curves and from unramified abelian covers of curves of genus $2$ and $3$. The largest example they gave was of genus $1297=1+6^4$, which was constructed as an unramified degree $1296$ abelian cover of a genus-two curve whose Jacobian is isogenous to a product of two elliptic curves \cite[\S4]{ES93}. Over the past three decades, much work has been done on trying to find completely decomposable Jacobians in genera not found by Ekedahl and Serre (called \textit{Ekedahl--Serre gaps}). Yamauchi \cite{Yam} gave a list of integers $N$ such that the Jacobian of the modular curve $X_0(N)$ is completely decomposable over $\mathbb{Q}$ and filled in several Ekedahl--Serre gaps. Paulhus and Rojas \cite{PR17} found many further genera using group actions and intermediate quotients. Rodr\'{\i}guez and Rojas \cite{RodRoj} used a general period matrix decomposition to fill in the gap $101$. Paulhus and Sutherland \cite{PS25} recently searched all modular curves $X_H$ of level $N<240$ for Jacobians that are isogenous \emph{over $\Q$} to products of elliptic curves, and found $13$ new genera, the largest being $409$. Other aspects of the problem, and its relation to other conjectures, have been discussed in \cite{ChenLuZuo} and \cite{ChenLuZuo2}; see also \cite[\S3.11]{Ser20} for a summary of the problem.

According to the list given in \cite{PS25}, the only known numbers $g$ over $300$ for which there exists a complex completely decomposable Jacobian of dimension $g$ are $325$, $409$, $433$, $649$, and $1297$. We also note that, despite the interest this problem has generated over the years, and the many contributions this area has received, Ekedahl and Serre's record of $1297$ remained the largest known genus until now.

\vspace{0.2cm}

Let $\h$ be the upper half-plane and define $\Gamma\subset\SL(2,\Z)$ to be the subgroup
\begin{equation}\label{eq:Gamma}
\Gamma:=\Bigl\{\gamma\in\SL(2,\Z)\;:\; \gamma\equiv\lambda I \ (\mathrm{mod}\ 8),\ \gamma\equiv\pm\mat1b01\ (\mathrm{mod}\ 3),\ \gamma\equiv\pm I\ (\mathrm{mod}\ 5)\Bigr\},
\end{equation}
where $\lambda$ and $b$ are not fixed and depend on $\gamma$.
The three conditions are imposed independently: the scalar $\lambda$ is necessarily a unit modulo $8$, and the signs modulo $3$ and $5$ are unrelated. Let $X_\Gamma$ be the compact Riemann surface obtained from $\Gamma\backslash\h$ by adding the cusps. Our main theorem shows that this curve breaks the Ekedahl--Serre record:

\begin{theorem}\label{maintheorem}
The curve $X_\Gamma$ has genus $1777$, admits a model over $\Q$, and $J(X_\Gamma)$ is completely decomposable over $\Qbar$.
\end{theorem}

We note that this decomposition is defined over $\Qbar$, and for the model of $X_\Gamma$ over $\Q$ used below it cannot be achieved over $\Q$ (Remark~\ref{rem:overQ}). This is also the reason the curve does not appear in the search of \cite{PS25}, which requires a splitting over $\Q$ into elliptic curves.

The idea of the proof is the following: We first find another model for $X_\Gamma$ that is more convenient to work with. We realize the new model as an \'etale double cover of the modular curve $X_0(4800)$. Using this double cover, we decompose, modulo isogeny, the Jacobian of $X_\Gamma$ into the product of the Jacobian of $X_0(4800)$ with another explicit abelian variety. Using the fact that, for a modular Jacobian, every $\Q$-simple factor can be described using certain weight-two newforms, we proceed to study these factors and show that they can all be decomposed even further into the product of elliptic curves over $\Qbar$. 

As an immediate corollary to the proof of Theorem \ref{maintheorem}, we obtain the following:

\begin{corollary}\label{cor:X0}
The modular curve $X_0(4800)$ has genus $889$, and its Jacobian is completely decomposable. This fills in another Ekedahl--Serre gap that was previously unknown.
\end{corollary}

Indeed, by Remark~\ref{rem:unramified} the curve $X_\Gamma$ is a double cover of $X_0(4800)$, so $J_0(4800)$ is an isogeny factor of $J(X_\Gamma)$; Corollary~\ref{cor:X0} is also proved directly in Section~\ref{sec:trivial}. It does not contradict \cite{Yam}, which determines the levels $N$ for which $J_0(N)$ is isogenous \emph{over $\Q$} to a product of elliptic curves: $J_0(4800)$ is not, since it has $\Q$-simple factors of dimension $2$ (Remark~\ref{rem:overQ}).

\subsection*{Notation}
For a smooth projective curve $X$ we write $J(X)$ for its Jacobian. An \emph{isogeny} is a surjective homomorphism of abelian varieties with finite kernel; we write $A\sim B$ if $A$ and $B$ are isogenous over an algebraically closed field, and $A\sim_k B$ if an isogeny exists over the field $k$. For an abelian variety $A$ over a field $k$, $\End(A)$ is the ring of endomorphisms of $A$ defined over $k$, and $\End_\Q(A):=\End(A)\otimes_\Z\Q$ is its endomorphism algebra; for an extension $k'/k$, $A_{k'}$ denotes the base change, so that $\End_\Q(A_{\Qbar})$ is the algebra of all endomorphisms of $A$ when $k$ is a number field. For a field $F$, $\Br(F)$ denotes the Brauer group of $F$. For $a,b\in\Q^\times$, $(a,b)_\Q$ is the quaternion algebra over $\Q$ with basis $1,u,v,uv$ and relations $u^2=a$, $v^2=b$, $uv=-vu$; it is \emph{split} if it is isomorphic to $M_2(\Q)$, which happens if and only if $a$ is a norm from $\Q(\sqrt b)$. For a discriminant $D$, $\chi_D=\left(\frac D\cdot\right)$ is the associated Kronecker character, and $\sigma_0(n)$ is the number of positive divisors of $n$. For $N\geq1$ we use the congruence subgroups
\begin{align*}
\Gamma_0(N)&=\Bigl\{\smat abcd\in\SL(2,\Z):c\equiv0\ (\mathrm{mod}\ N)\Bigr\},\\
\Gamma_1(N)&=\Bigl\{\smat abcd\in\Gamma_0(N):a\equiv d\equiv1\ (\mathrm{mod}\ N)\Bigr\},
\end{align*}
the modular curves $X_0(N)$ and $X_1(N)$, obtained from $\Gamma_0(N)\backslash\h$ and $\Gamma_1(N)\backslash\h$ by adding the cusps, and their Jacobians $J_0(N)=J(X_0(N))$ and $J_1(N)=J(X_1(N))$. We write $\h^*=\h\cup\PP^1(\Q)$.

\subsection*{Acknowledgements}
The author was partially supported by the University of Chile Enlace Grant ENL31/26.

\subsection*{Use of artificial intelligence}
AI tools were used in exploratory research and in writing and cross-checking the computer verifications, as well as in improving the exposition of the results. The author is responsible for the writing and the accuracy of the results in the article.

\section{The curve $X_\Gamma$}

In this section we will prove basic properties of the curve $X_\Gamma$. For $L\geq 1$, let $\Gamma(L)$ be the kernel of the quotient map $\SL(2,\Z)\to\SL(2,\Z/L\Z)$. Since the definition of the elements in $\Gamma$ depends only on $\gamma$ modulo $120=3\cdot5\cdot 8$, we have that $\Gamma(120)\subseteq\Gamma$. By the Chinese Remainder Theorem, we have that
\[\SL(2,\Z/120\Z)\simeq\SL(2,\Z/8\Z)\times\SL(2,\Z/3\Z)\times\SL(2,\Z/5\Z),\]
and the image $G_\Gamma$ of $\Gamma$ in $\SL(2,\Z/120\Z)$ is isomorphic to the product of the following three groups:
\[G_3:=\left\{\pm\begin{pmatrix}1&b\\0&1\end{pmatrix}:b\in\Z/3\Z\right\}\]
\[G_5:=\{\pm I\}\]
\[G_8:=\{\lambda I:\lambda\in(\Z/8\Z)^\times\}.\]
We note that $-I\in\Gamma$. Now if $G$ is any finite-index subgroup of $\SL(2,\Z)$ containing $-I$ of index $\mu$ with $e_p$ elliptic points of order $p$ and $c$ cusps, the Riemann--Hurwitz formula for $X_G\to X(1)$ gives us (see \cite[Theorem 3.1.1]{DS05}):
\begin{equation}\label{genus} g(X_G)=1+\frac{\mu}{12}-\frac{e_2}{4}-\frac{e_3}{3}-\frac{c}{2}.\end{equation}

We then have the following:

\begin{proposition}
The group $\Gamma$ has level $120$, index $\mu=23040$, no elliptic points, and $288$ cusps. In particular, $g(X_\Gamma)=1777$.
\end{proposition}
\begin{proof}
Given that $|\SL(2,\Z/p^e\Z)|=p^{3e-2}(p^2-1)$ for any prime $p$ and positive integer $e$, we have that 
\[|\SL(2,\Z/120\Z)|=1105920\]
\[|G_\Gamma|=48.\]
In particular, $\mu=1105920/48=23040$. Now we see that if $\Gamma(L)\subseteq\Gamma$, then in particular $\smat 10L1\in\Gamma$, which implies that $8\mid L$, $3\mid L$ and $5\mid L$. In particular, $120\mid L$, and therefore we conclude that the level of $\Gamma$ is exactly $120$.
As for the elliptic elements of $\Gamma$, we have that an elliptic element of $\SL(2,\Z)$ must have trace equal to either $0$ or $\pm1$. However, since every element of $\Gamma$ is congruent to a scalar matrix $\lambda I$ modulo $8$ such that $\lambda$ is odd, we have that its trace must be even and non-zero. This implies that $\Gamma$ contains no elliptic elements, and so $e_2=e_3=0$. 

Finally, we have that the cusps of $\Gamma$ are associated with the $\Gamma$-orbits on primitive vectors in $\Z^2$ up to sign. Since $\Gamma$ contains $\Gamma(120)$ and $-I$, we have that these orbits therefore correspond to $G_\Gamma$-orbits on the set of primitive vectors in $(\Z/120\Z)^{2}$. This implies that we are looking for vectors in $(\Z/120\Z)^2$ whose coordinates generate the unit ideal. We can look at these sets modulo $3$, $5$ and $8$: Modulo 3, on the 8 non-zero vectors $(x,y)$, the unipotent matrices change $x$ freely when $y\neq0$, and the sign identifies $y$ with $-y$. This gives us one orbit with $y\neq 0$ and one with $y=0$, and therefore a total of $2$ orbits. Modulo $5$, the $24$ non-zero vectors lie in exactly $12$ orbits under $\pm1$. Modulo 8, there are $8^2-4^2=48$ primitive vectors and the four scalars act freely, which gives us exactly $12$ orbits. Therefore, the number of total orbits, and therefore cusps, is exactly $c=2\cdot 12\cdot 12=288$. The result follows by substituting these values into \eqref{genus}.
\end{proof}

If $N\geq 1$ and $H\leq(\Z/N\Z)^\times$ is a subgroup containing $-1$, define
\[\Gamma_H(N):=\left\{\begin{pmatrix}a&b\\c&d\end{pmatrix}\in\Gamma_0(N):d\text{ mod }N\in H\right\}.\]
We see that if $H=(\Z/N\Z)^\times$, then $\Gamma_H(N)=\Gamma_0(N)$, and $\Gamma_{\{\pm1\}}(N)=\pm\Gamma_1(N)$. 

Now set $N=4800=2^6\cdot3\cdot5^2$ and $H=\{d\in(\Z/4800\Z)^\times:d\equiv\pm1\pmod{5}\}$. 

\begin{proposition}\label{prop:conjugate}
Consider the element $\alpha=\smat{40}001\in\GL^+(2,\Q)$. Then $\alpha^{-1}\Gamma\alpha=\Gamma_H(4800)$, and so the map $z\mapsto z/40$ induces an isomorphism of Riemann surfaces $X_\Gamma\to X_{\Gamma_H(4800)}$.
\end{proposition}
\begin{proof}
For $\gamma=\smat abcd$ we have $\alpha^{-1}\gamma\alpha=\smat{a}{b/40}{40c}{d}$. If $\gamma\in\Gamma$, then modulo $8$ and modulo $5$, $\gamma$ is scalar, so $b\equiv0\pmod{40}$ and $b/40\in\Z$. Moreover $v_2(c)\geq3$, $v_5(c)\geq1$ and, since $\gamma$ is upper triangular modulo $3$, $3\mid c$; hence $40c\equiv0\pmod{2^6\cdot3\cdot5^2}$. Finally $d\equiv\pm1\pmod5$, which implies that $\alpha^{-1}\gamma\alpha\in\Gamma_H(4800)$.
Conversely, let $\gamma'=\smat abcd\in\Gamma_H(4800)$ and put $\gamma=\alpha\gamma'\alpha^{-1}=\smat{a}{40b}{c/40}{d}$. Since $4800\mid c$, we have $120\mid c/40$, so $\gamma\in\SL(2,\Z)$ is upper triangular modulo $120$ with $ad\equiv1\pmod{120}$. Modulo $8$, every unit is its own inverse, so $a\equiv d$ and $\gamma\equiv aI$ is scalar. Modulo $3$, $a\equiv d\equiv\pm1$, so $\gamma\equiv\pm\smat1*01$. Modulo $5$, $d\equiv\pm1$ and $a\equiv d^{-1}\equiv d$, while $40b\equiv0$; so $\gamma\equiv\pm I$. Hence $\gamma\in\Gamma$.
\end{proof}

\begin{remark}\label{rem:unramified}
Since $[\Gamma_0(4800):\Gamma_H(4800)]=[(\Z/4800\Z)^\times:H]=2$, the curve $X_\Gamma$ is a double cover of $X_0(4800)$. The genus of $X_0(4800)$ is $889$ and it has $144$ cusps \cite[\S3.8]{DS05}, so Riemann--Hurwitz shows that the cover is unramified; in particular each cusp of $X_0(4800)$ splits into two of the $288$ cusps of $X_\Gamma$.
\end{remark}

The modular curve $X_1(N)$ has a standard smooth projective model over $\Q$, on which the diamond operators $\diam d$, $d\in(\Z/N\Z)^\times$, act by automorphisms defined over $\Q$ \cite[Chapter~7]{DS05}. We define
\[
X_H(N):=X_1(N)/\{\diam d:d\in H\},
\]
a smooth projective curve over $\Q$ with $X_H(N)(\C)\cong\Gamma_H(N)\backslash\h^*$. By Proposition~\ref{prop:conjugate}, $X_H(4800)$ is a model of $X_\Gamma$ over $\Q$, and we use this model throughout. Since Theorem~\ref{maintheorem} concerns $J(X_\Gamma)$ over $\Qbar$, the choice of model plays no role in the decomposition.

\begin{remark}\label{rem2.4}
Following \cite[\S2 and Appendix~A]{RSZB}, $X_\Gamma$ also has the model $X_{H_{120}}$ over $\Q$ attached to the subgroup $H_{120}=\langle G_\Gamma,\ \smat u001:u\in(\Z/120\Z)^\times\rangle\subseteq\GL(2,\Z/120\Z)$, which has surjective determinant and satisfies $H_{120}\cap\SL(2,\Z/120\Z)=G_\Gamma$. The Frobenius traces of $J(X_{H_{120}})$, computed by an independent point count, agree with those of $J(X_H(4800))$ at all primes $7\leq p\leq43$ (Table~\ref{tab:traces}). We do not need to compare the two models further.
\end{remark}

\section{Preliminaries on modular abelian varieties}\label{sec:background}

In this section we collect standard results on modular abelian varieties that we will use in later sections. 

If $G$ is any congruence subgroup, the map $f\mapsto f(z)dz$ identifies the space of weight-two cusp forms $S_2(G)$ with $H^0(X_G,\Omega^1)$, which is just the cotangent space of $J(X_G)$ at the origin. On $X_1(N)$, the Hecke correspondences $T_p$ ($p\nmid N$) and the diamond automorphisms $\langle d\rangle$ are defined over $\Q$, and act compatibly on the Jacobian and on the space $S_2(\Gamma_1(N))$.

If $\varepsilon$ is a Dirichlet character modulo $M$ with $\varepsilon(-1)=1$, let $S_2(M,\varepsilon)$ denote the space of weight 2 cusp forms on $\Gamma_1(M)$ with $f|\langle d\rangle=\varepsilon(d)f$, and let $S_2^{\operatorname{new}}(M,\varepsilon)$ be its new subspace. Let $q=e^{2\pi i z}$ be the standard parameter for $z\in\h$, and for a prime $p$ define the operator
\[U_p\Bigl(\sum b_n q^n\Bigr):=\sum b_{pn}q^n.\]
A \textit{newform} is a normalized ($a_1=1$) eigenform $f=\sum a_nq^n\in S_2^{\operatorname{new}}(M,\varepsilon)$ for all Hecke operators, which implies that $T_pf=a_pf$ for every $p$ that doesn't divide $M$ and $U_pf=a_pf$ for $p$ dividing $M$. Let $K_f:=\Q(a_n:n\geq1)$ denote the coefficient field of $f$. We have that $\operatorname{Gal}(\Qbar/\Q)$ permutes newforms by acting on their coefficients; denote by $[f]$ the orbit of $f$ under the action of the Galois group. 

We have the following properties.

\subsection{Strong multiplicity one}\label{s1} (see \cite{AL70,Li75}): A newform is determined by its eigenvalues $a_p$ at all but finitely many primes. This is called the \textit{strong multiplicity one} property for newforms, and is equivalent to saying that in $S_2(L,\varepsilon)$, the forms with prescribed eigenvalues for all $T_p$, $p\nmid L$, form a space that contains exactly one newform. If $M$ denotes the level of this newform, then $M$ divides $L$, and the space is spanned by $f(q^e)$ for $e\mid L/M$. If $\varepsilon=1$ and $p\mid M$, then $a_p\neq0$ if and only if the largest power of $p$ that divides $M$ is $p$. Consequently,
\[\dim S_2(\Gamma_1(N))=\sum\sigma_0(N/M),\]
where the sum ranges over all newforms of level $M\mid N$.

\subsection{The abelian variety attached to a newform}\label{s2} (see \cite[Theorem 7.14]{Shi71}): To every orbit $[f]$, where $f$ is a level $M$ newform, we can attach an abelian variety $A_f$ over $\Q$ which is a quotient of the Jacobian $J_1(M)$ of the modular curve $X_1(M)$. We have that
\begin{itemize}
\item $\dim A_f=[K_f:\Q]$
\item $\End_\Q(A_f)\simeq K_f$; in particular $A_f$ is simple over $\Q$. 
\item $A_f$ is not necessarily geometrically simple, and there exists a simple abelian variety $B_f$ such that, over $\Qbar$, $A_f$ is isogenous to a power of $B_f$. In particular, $A_f$ is completely decomposable if and only if $B_f$ is an elliptic curve. We call $B_f$ the \textit{building block} of $f$ (see \cite{Rib80,Rib92} and Theorem~\ref{theoremresults} below).
\item For every prime $\ell$,
\begin{equation}\label{Vell} V_\ell(A_f)\otimes_{\Q_\ell}\Qbar_\ell\simeq\bigoplus_{\tau:K_f\hookrightarrow\Qbar_\ell}\rho_{f,\tau},\end{equation}
where $\rho_{f,\tau}$ is the semisimple two-dimensional Galois representation with $\operatorname{Tr}\rho_{f,\tau}(\operatorname{Frob}_p)=\tau(a_p)$ for every $p\neq\ell$ that does not divide $M$.
\end{itemize}

\subsection{Twists of newforms}\label{s3} (see \cite[\S3]{AL78}): If $\chi$ is a primitive Dirichlet character of conductor $r$, then we define the \textit{raw twist} of a newform $f=\sum a_nq^n\in S_2(M,\varepsilon)$ to be 
\[f_\chi:=\sum \chi(n)a_nq^n\in S_2\bigl(\lcm(M,r^2,r\cdot\cond(\varepsilon)),\varepsilon\chi^2\bigr).\]
For every prime $p$ that does not divide $Mr$, it has eigenvalue $\chi(p)a_p$ for $T_p$. By \ref{s1}, the space of forms with these eigenvalues contains exactly one newform which is called the \textit{primitive twist} of $f$ and is denoted by $f\otimes\chi$. The forms $f_\chi$ and $f\otimes\chi$ have the same eigenvalues at the primes $p\nmid Mr$, but in general they differ at the primes dividing $Mr$. We will be using raw twists when we need explicit $q$-expansions, and primitive twists when we need to use newforms.

\subsection{Bound on the coefficients}\label{s4} (see \cite{Stu87} and \cite[\S8.2]{BBB+}): If we define
\[B(L):=\left\lceil \frac{1}{6}[\SL(2,\Z):\Gamma_0(L)]\right\rceil,\]
then two forms in $S_2(L,\varepsilon)$ are equal if and only if for all $n\leq B(L)$, their $n$th coefficients agree.

We now proceed to decompose the Jacobian of the modular curve $X_H(N)$:

\begin{theorem}\label{decomposition}
Let $H\leq (\Z/N\Z)^\times$ be a subgroup that contains $-1$. If $f$ is a newform of level $M\mid N$ and $\varepsilon_f$ is its character (regarded as a character modulo $N$), then
\[J(X_H(N))\sim\prod_{M\mid N}\prod_{[f]}A_f^{\sigma_0(N/M)}.\]
Here, for each $M$, the inner product runs over the $\text{Gal}(\Qbar/\Q)$-orbits of newforms $[f]$ of level $M$ such that $\varepsilon_f$ is trivial on $H$, and the isogeny is defined over $\Q$.
\end{theorem}
\begin{proof}
If $H=\{\pm1\}$ then this is exactly the content of \cite[Theorem 6.6.6]{DS05}. In general, we have that $X_H(N)=X_1(N)/\tilde{H}$, where $\tilde{H}$ is the image of $H$ in the automorphism group of $X_1(N)$. By pulling back, we have that up to isogeny over $\Q$, $J(X_H(N))$ can be identified with the image of $e=\sum_{d\in H}\langle d\rangle$. On $A_f^{\sigma_0(N/M)}$, the operator $\langle d\rangle$ acts through $\varepsilon_f(d)\in K_f$, and so $e$ acts as multiplication by $|H|$ if $\varepsilon_f$ is trivial on $H$ and $0$ if not.
\end{proof}

We will now look at the building blocks of $A_f$. A newform $f$ has CM (complex multiplication) by $\Q(\sqrt{D})$, for $D<0$, if $a_p=\chi_D(p)a_p$ for almost all $p$, where $\chi_D=\left(\frac{D}{\cdot}\right)$. An \emph{inner twist} of $f$ is a pair $(\sigma,\chi)$, where $\sigma$ is an automorphism of $K_f$ and $\chi$ is a Dirichlet character that satisfies that $\sigma(f)=f\otimes\chi$. We recall the following theorem which brings together results proved by Shimura, Ribet, Momose, Quer, among others:

\begin{theorem}\label{theoremresults}
Let $f$ be a newform of weight $2$. Then we have the following:
\begin{enumerate}
\item If $f$ has CM, then $(A_f)_{\Qbar}\sim E^{[K_f:\Q]}$, where $E$ is an elliptic curve with CM (see \cite{Shi71b,Rib80}).
\item If $f$ has no CM and $F$ is the subfield of $K_f$ fixed by the automorphisms that admit an inner twist, then $(A_f)_{\Qbar}\sim B_f^r$ for a simple abelian variety $B_f$ such that $\dim B_f=t[F:\Q]$ for $t\in\{1,2\}$. Moreover, $t=1$ if and only if the class of the central simple $F$-algebra $\End_\Q((A_f)_{\Qbar})$ in $\Br(F)$ is trivial (see \cite{Rib80,Mom81,Rib92,Pyl04}).
\item If $f$ does not have CM, $\varepsilon=1$, $K_f=\Q(\sqrt{m})$ is real quadratic and $(\sigma,\chi_D)$ is an inner twist with non-trivial $\sigma$, then $F=\Q$, and the class of $\End_\Q((A_f)_{\Qbar})$ in $\Br(\Q)$ is that of the quaternion algebra $(m,D)_\Q$ (see \cite{Que09}).
\end{enumerate}
\end{theorem}

\begin{remark}\label{rem:quaternion} We note that, in particular, if $K_f$ is a real quadratic field $\Q(\sqrt{m})$, $\varepsilon=1$, $f$ has no CM and $(\sigma,\chi_D)$ is an inner twist with $\sigma\neq\text{id}$, then $B_f$ is an elliptic curve if and only if $(m,D)_\Q$ is split. 
\end{remark}

We will need the following two lemmas in the sections that follow:

\begin{lemma}\label{lemma1}
Let $f$ and $g$ be newforms of weight $2$, and let $\chi$ be a Dirichlet character that satisfies $a_p(g)=\chi(p)a_p(f)$ for almost every $p$. Then $(A_f)_{\Qbar}$ and $(A_g)_{\Qbar}$ have the same simple factors up to isogeny (although the multiplicities could differ). In particular, $B_f\sim B_g$.
\end{lemma}
\begin{proof}
Let $K_\chi$ be the abelian field cut out by $\chi$, and let $L\supseteq K_\chi$ be a finite extension contained in $\Qbar$. By the Chebotarev Density Theorem and the Brauer-Nesbitt Theorem, we have that $\rho_{g,\tau'}\simeq\rho_{f,\tau}\otimes\chi$ for compatible embeddings $\tau$ and $\tau'$, which implies that $\rho_{g,\tau'}$ and $\rho_{f,\tau}$ agree on $\text{Gal}(\Qbar/L)$. By Equation \ref{Vell}, $V_\ell(A_f)$ and $V_\ell(A_g)$ have the same irreducible constituents as $\Qbar_\ell[\text{Gal}(\Qbar/L)]$-modules. If $B$ is a simple abelian subvariety of $(A_f)_{\Qbar}$ and defined over $L$, then $\text{Hom}_{\text{Gal}(\Qbar/L)}(V_\ell B,V_\ell A_g)\neq 0$. By Faltings' theorem \cite[Corollary~1]{Fal83}, there exists a non-zero homomorphism $B\to (A_g)_L$ whose image is isogenous to $B$. The claim follows by symmetry.
\end{proof}

\begin{lemma}\label{lemma2}
Let $f=\sum a_nq^n\in S_2^{\operatorname{new}}(M,1)$ be a newform, let $\sigma$ be a field automorphism of $K_f$, and let $D$ be a fundamental discriminant whose prime divisors all divide $M$. If $L=\operatorname{lcm}(M,D^2)$, and $\sigma(a_n)=\chi_D(n)a_n$ for every $n\leq B(L)$ with $\operatorname{gcd}(n,D)=1$, then $\sigma(a_p)=\chi_D(p)a_p$ for every prime $p$ that does not divide $D$. In particular, $(\sigma,\chi_D)$ is an inner twist. If $\sigma$ is the identity and $D<0$, then $f$ has CM by $\Q(\sqrt{D})$. 
\end{lemma}
\begin{proof}
Define the operators
\[V_p\left(\sum b_nq^n\right):=\sum b_n q^{pn}\]
\[\mathcal{D}_p\left(\sum b_nq^n\right)=\sum_{p\nmid n}b_nq^n.\]
Then we see that $\mathcal{D}_p=1-V_pU_p$. For the newform $\sigma(f)$, if $p\mid D$, then $\mathcal{D}_p\sigma(f)=\sigma(f)-\sigma(a_p)V_p(\sigma(f))$. By \ref{s1}, this second term appears only when the highest power of $p$ that divides $M$ is $1$. We have that $g_1=\prod_{p\mid D}\mathcal{D}_p\sigma(f)$ is a combination of the forms $V_e\sigma(f)$ with $e$ a product of primes dividing $D$ and that appear in $M$ with exponent $1$. This implies that $g_1\in S_2(L,1)$, and, by \ref{s3}, the raw twist $g_2=f_{\chi_D}$ also belongs to $S_2(L,1)$. If the greatest common divisor of $n$ and $D$ is $1$, then the $n$th coefficient of $g_1$ and $g_2$ are $\sigma(a_n)$ and $\chi_D(n)a_n$, respectively, and vanish otherwise. By \ref{s4}, we obtain that $g_1=g_2$, and the lemma follows by comparing coefficients at primes that do not divide $D$.
\end{proof}

\section{Analyzing the modular constituents}\label{sec:reduction}

We now apply Theorem~\ref{decomposition} to $X_H(4800)$ and remove the character $\chi_5$ by a quartic twist.

The characters that are trivial on $H=\{d:d\equiv\pm1\bmod5\}$ modulo $4800$ are those that factor through the group $(\Z/4800\Z)^\times/H\cong(\Z/5\Z)^\times/\{\pm1\}$. There are exactly two: the trivial character and $\chi_5=(\tfrac5\cdot)$. Theorem~\ref{decomposition} therefore gives us:

\[
J(X_H(4800))\sim\prod_{M\mid4800}\ \prod_{[f]\subset S_2^\new(M,1)}A_f^{\sigma_0(4800/M)}\times\prod_{5\mid M\mid4800}\ \prod_{[f]\subset S_2^\new(M,\chi_5)}A_f^{\sigma_0(4800/M)},\]
where the first double product corresponds to $J_0(4800)$, the Jacobian of the modular curve $X_0(4800)$, and the second product, which we denote by $J_{\chi_5}$, is an abelian variety of dimension $888$ over $\Q$. On cotangent spaces this corresponds precisely to the decomposition
\[S_2(\Gamma_H(4800))=S_2(4800,1)\oplus S_2(4800,\chi_5).\] 
The two summands have dimensions $889$ and $888$, respectively. Note that the first factor is isogenous to $J_0(4800)$ by the same theorem with the group $H=(\Z/4800\Z)^\times.$ The dimensions are computed from the standard dimension formulas, or directly from the data in Appendix~\ref{app:spaces}. 

Now let $\psi$ be a Dirichlet character modulo $5$ of order $4$. Then $\psi^2=\chi_5$, $\psi$ is primitive of conductor $5$, and $\psi(-1)=-1$.

\begin{proposition}\label{prop:quartic}
Let $f$ be a newform of level $M\mid4800$ and character $\chi_5$. Then the primitive twist $g=f\otimes\psi^{-1}$ is a newform of trivial character and level dividing $4800$, that satisfies $a_p(g)=\psi^{-1}(p)\,a_p(f)$ for all primes $p$ that do not divide $5M$. In particular, $B_f\sim B_g$.
\end{proposition}

\begin{proof}
Here we see that $r=\cond(\psi)=5$ and $\cond(\chi_5)=5$, which implies that the raw twist satisfies
\[
f_{\psi^{-1}}\in S_2\bigl(\lcm(M,25),\chi_5\psi^{-2}\bigr)=S_2\bigl(\lcm(M,25),1\bigr),
\]
and it is an eigenform of $T_p$ for every $p\nmid5M$ with eigenvalues $\psi^{-1}(p)a_p(f)$. This implies that the primitive twist $g=f\otimes\psi^{-1}$ is a newform of trivial character whose level divides $\lcm(M,25)$, which divides $4800$. The raw twist itself is not necessarily a newform. The last assertion is Lemma~\ref{lemma1}.
\end{proof}

\begin{corollary}\label{cor:reduction}
Every geometric simple factor of $J_{\chi_5}$ is isogenous to a geometric simple factor of $J_0(4800)$. In particular, if $J_0(4800)$ is completely decomposable over $\Qbar$, then so is $J(X_\Gamma)$.
\end{corollary}

\begin{proof}
Combine the decomposition of the Jacobian $J(X_H(4800))$ with Proposition \ref{prop:quartic}.
\end{proof}

\section{Newforms of trivial character and level dividing 4800}\label{sec:trivial}

By Corollary~\ref{cor:reduction}, we see that Theorem~\ref{maintheorem} is an immediate consequence of the complete decomposition of $J_0(4800)$; we will now proceed with the proof. This is the only step that requires computation, and the computation has two parts:
\begin{enumerate}[label=(C\arabic*)]
\item the decomposition of the newspaces $S_2^\new(M,1)$, $M\mid4800$, into Galois orbits, together with their coefficient fields;
\item for each of the non-rational orbits, the comparison of $q$-expansions through a Sturm bound required by Lemma~\ref{lemma2}.
\end{enumerate}
Both are exact computations; they are described, together with independent cross-checks, in Appendix~\ref{app:repro}.

Decomposing the $30$ nonzero newspaces $S_2^\new(M,1)$, where $M$ is a divisor of $4800$, into Galois orbits (computation (C1); Appendix~\ref{app:spaces}) gives $264$ orbits. Of these:
\begin{itemize}
\item $248$ have $K_f=\Q$;
\item $4$ have $K_f=\Q(\sqrt2)$, at levels $160$, $320$, $800$ and $1600$;
\item $12$ have $K_f=\Q(\sqrt5)$, at levels $800$ (four orbits), $1600$ (four), $2400$ (two) and $4800$ (two).
\end{itemize}
No other coefficient fields occur. Table~\ref{tab:sixteen} lists the sixteen non-rational orbits. Each is determined by its level and by the eigenvalues $a_7,a_{11},a_{13},a_{17}$ given there, up to $\sqrt m\mapsto-\sqrt m$.

\begin{table}[ht]
\centering\footnotesize
\setlength{\tabcolsep}{2.5pt}
\caption{The sixteen Galois orbits of newforms of weight $2$, trivial character and level $M\mid4800$ with $K_f\neq\Q$. ``Type'' gives the self-twist (CM) or inner twist proved by Lemma~\ref{lemma2}, with $L=\lcm(M,D^2)$ and Sturm bound $B(L)$. The last column is the multiplicity $\sigma_0(4800/M)$ of $A_f$ in $J_0(4800)$.}\label{tab:sixteen}
\begin{tabular}{rcccccccccc}
\toprule
$M$ & $K_f$ & $a_7$ & $a_{11}$ & $a_{13}$ & $a_{17}$ & type & $L$ & $B(L)$ & symbol & mult.\\
\midrule
160 & $\Q(\sqrt{2})$ & $2\sqrt{2}$ & $4\sqrt{2}$ & $-2$ & $2$ & twist $\chi_{-4}$ & 160 & 48 & $(2,-4)$ & 8 \\
320 & $\Q(\sqrt{2})$ & $-2\sqrt{2}$ & $4\sqrt{2}$ & $2$ & $2$ & twist $\chi_{-4}$ & 320 & 96 & $(2,-4)$ & 4 \\
800 & $\Q(\sqrt{5})$ & $3+\sqrt{5}$ & $0$ & $0$ & $0$ & CM, $D=-20$ & 800 & 240 & -- & 4 \\
800 & $\Q(\sqrt{2})$ & $2\sqrt{2}$ & $-4\sqrt{2}$ & $2$ & $-2$ & twist $\chi_{-4}$ & 800 & 240 & $(2,-4)$ & 4 \\
800 & $\Q(\sqrt{5})$ & $-3+\sqrt{5}$ & $0$ & $0$ & $0$ & CM, $D=-20$ & 800 & 240 & -- & 4 \\
800 & $\Q(\sqrt{5})$ & $-2\sqrt{5}$ & $-\sqrt{5}$ & $4$ & $-7$ & twist $\chi_{-4}$ & 800 & 240 & $(5,-4)$ & 4 \\
800 & $\Q(\sqrt{5})$ & $-2\sqrt{5}$ & $\sqrt{5}$ & $-4$ & $7$ & twist $\chi_{-4}$ & 800 & 240 & $(5,-4)$ & 4 \\
1600 & $\Q(\sqrt{5})$ & $-3-\sqrt{5}$ & $0$ & $0$ & $0$ & CM, $D=-20$ & 1600 & 480 & -- & 2 \\
1600 & $\Q(\sqrt{2})$ & $-2\sqrt{2}$ & $-4\sqrt{2}$ & $-2$ & $-2$ & twist $\chi_{-4}$ & 1600 & 480 & $(2,-4)$ & 2 \\
1600 & $\Q(\sqrt{5})$ & $3-\sqrt{5}$ & $0$ & $0$ & $0$ & CM, $D=-20$ & 1600 & 480 & -- & 2 \\
1600 & $\Q(\sqrt{5})$ & $-2\sqrt{5}$ & $-\sqrt{5}$ & $4$ & $7$ & twist $\chi_{-4}$ & 1600 & 480 & $(5,-4)$ & 2 \\
1600 & $\Q(\sqrt{5})$ & $-2\sqrt{5}$ & $\sqrt{5}$ & $-4$ & $-7$ & twist $\chi_{-4}$ & 1600 & 480 & $(5,-4)$ & 2 \\
2400 & $\Q(\sqrt{5})$ & $2$ & $-2\sqrt{5}$ & $-2\sqrt{5}$ & $2\sqrt{5}$ & twist $\chi_{-20}$ & 2400 & 960 & $(5,-20)$ & 2 \\
2400 & $\Q(\sqrt{5})$ & $-2$ & $-2\sqrt{5}$ & $2\sqrt{5}$ & $-2\sqrt{5}$ & twist $\chi_{-20}$ & 2400 & 960 & $(5,-20)$ & 2 \\
4800 & $\Q(\sqrt{5})$ & $-2$ & $-2\sqrt{5}$ & $2\sqrt{5}$ & $2\sqrt{5}$ & twist $\chi_{-20}$ & 4800 & 1920 & $(5,-20)$ & 1 \\
4800 & $\Q(\sqrt{5})$ & $2$ & $-2\sqrt{5}$ & $-2\sqrt{5}$ & $-2\sqrt{5}$ & twist $\chi_{-20}$ & 4800 & 1920 & $(5,-20)$ & 1 \\
\bottomrule
\end{tabular}
\end{table}

\subsection{The CM orbits}
For four of the $\Q(\sqrt5)$-orbits, at levels $800$ and $1600$, the coefficients $a_p$ vanish at every good prime $p<400$ with $\chi_{-20}(p)=-1$. Lemma~\ref{lemma2} with $\sigma=\mathrm{id}$ and $D=-20$, checked through $B(800)=240$ and $B(1600)=480$ coefficients, proves that these forms have CM by $\Q(\sqrt{-5})$. By Theorem~\ref{theoremresults}(1), $(A_f)_{\Qbar}$ is isogenous to the square of an elliptic curve with CM by an order in $\Q(\sqrt{-5})$.

\subsection{The inner twists}
For each of the other twelve orbits, let $\sigma$ be the non-trivial automorphism of $K_f=\Q(\sqrt m)$, $m\in\{2,5\}$. Among the quadratic characters $\chi_D$ with $D$ supported on $\{2,3,5\}$, exactly one satisfies $\sigma(a_p)=\chi_D(p)a_p$ at all good primes $p<400$: $D=-4$ for the four $\Q(\sqrt2)$-orbits and for the four $\Q(\sqrt5)$-orbits of levels $800$ and $1600$, and $D=-20$ for the four $\Q(\sqrt5)$-orbits of levels $2400$ and $4800$. Lemma~\ref{lemma2} proves each of these inner twists (computation (C2)), using the Sturm bounds in Table~\ref{tab:sixteen}; the largest is $B(4800)=1920$. The pattern is visible in the table: $a_p\in\Q\sqrt m$ exactly when $\chi_D(p)=-1$.

Since $\sigma$ admits an inner twist, $F=\Q$. If $f$ had CM, it would have an elliptic building block by Theorem~\ref{theoremresults}(1); so we may assume it has no CM, and Theorem~\ref{theoremresults}(3) can be applied. The three symbols that occur are split:
\begin{align*}
(2,-4)&=(2,-1),\qquad &2&=N_{\Q(i)/\Q}(1+i),\\
(5,-4)&=(5,-1),\qquad &5&=N_{\Q(i)/\Q}(2+i),\\
(5,-20)&=(5,-5),\qquad &-5&=N_{\Q(\sqrt5)/\Q}(\sqrt5).
\end{align*}
Here we used that $(a,b)$ depends only on the square classes of $a$ and $b$, and that $(a,b)$ is split if and only if $a$ is a norm from $\Q(\sqrt b)$. Hence $\End_\Q((A_f)_{\Qbar})\cong M_2(\Q)$ and $(A_f)_{\Qbar}\sim E^2$ for an elliptic curve $E$.

\begin{proposition}\label{prop:allelliptic}
Every newform of weight $2$, trivial character and level dividing $4800$ has an elliptic building block.
\end{proposition}

\begin{proof}
For the $248$ rational orbits, $A_f$ is an elliptic curve over $\Q$. The sixteen remaining orbits are treated in the two preceding subsections.
\end{proof}

In particular, $J_0(4800)$ is completely decomposable, and the proof of Theorem \ref{maintheorem} is complete.

\subsection{Dimension count}
Table~\ref{tab:count} records how the $1777$ dimensions are distributed among the constituents of the Jacobian $J(X_\Gamma)$. Each orbit $[f]$ of level $M$ contributes $\sigma_0(4800/M)[K_f:\Q]$ dimensions, and each such contribution consists of that many elliptic factors over $\Qbar$.

\begin{table}[ht]
\centering\small
\caption{Distribution of the $1777$ dimensions of $J(X_\Gamma)$ by character and coefficient field. Here $\zeta_8$ is a primitive eighth root of unity, so $\Q(\zeta_8)=\Q(i,\sqrt2)$.}\label{tab:count}
\begin{tabular}{llrr}
\toprule
character & coefficient field & Galois orbits & dimension\\
\midrule
$1$ & $\Q$ & $248$ & $793$\\
$1$ & $\Q(\sqrt2)$ & $4$ & $36$\\
$1$ & $\Q(\sqrt5)$ & $12$ & $60$\\
\midrule
$\chi_5$ & $\Q(i)$ & $124$ & $768$\\
$\chi_5$ & $\Q(i,\sqrt5)$ & $6$ & $96$\\
$\chi_5$ & $\Q(\zeta_8)$ & $2$ & $24$\\
\midrule
total & & $396$ & $1777$\\
\bottomrule
\end{tabular}
\end{table}

\begin{remark}\label{rem:overQ}
Over $\Q$, the decomposition of $J_0(4800)$ contains abelian varieties $A_f$ with $K_f\neq\Q$, which are simple over $\Q$ because $\End_\Q(A_f)=K_f$ is a field. Thus $J(X_\Gamma)$ is not $\Q$-isogenous to a product of elliptic curves over $\Q$, for the model $X_H(4800)$; its largest $\Q$-split part is the product of the elliptic curves of the $248$ rational orbits, of dimension $793$. This is why the curve is not found by searches for $\Q$-split Jacobians such as \cite{PS25}. Allowing geometric splitting, as in \cite{ES93}, makes the twisting argument of Proposition~\ref{prop:quartic} available.
\end{remark}

\subsection{An independent verification}\label{sec:equivariant}
The decomposition was first obtained by a different, purely finite method. This first method does not use Proposition~\ref{prop:conjugate}, and it applies to any subgroup containing a principal congruence subgroup. We only sketch it here; the certificate and its verifiers are part of the supplementary material (Appendix~\ref{app:repro}).

One views $X_\Gamma$ as a quotient of the full modular curve of level $N_0=360$ (with all its $96$ geometric components), on which $G=\GL(2,\Z/360\Z)$ acts, and decomposes its space of holomorphic differentials into $G$-isotypic parts using the Chevalley--Weil formula \cite{CW34}. Grouping the representations that differ by a character of $\det$ gives $44$ blocks, whose contributions to the genus of $X_\Gamma$ add up to $1777$. For a representative of each block, explicit local invariant vectors and Casselman's newvector theorem \cite{Cas73} bound the conductors of all possible constituents. This leaves $38$ candidate newspaces of level dividing $2400$, which agree with the corresponding rows of Table~\ref{tab:spaces}. The Hecke traces of each representative are computed exactly in $\Q(\zeta_{24})$ from an equivariant point count over the Frobenius matrices of elliptic curves over $\F_p$ \cite{DT02}, for the $354$ primes $7\leq p\leq2399$. The multiplicities of the candidate newforms are then determined by exact rank certificates, namely nonzero maximal minors modulo $10^9+9$. The constituents of the block representatives are $49$ Galois orbits of newforms of level dividing $2400$ and character $1$ or $\chi_5$, all with elliptic building blocks; the remaining constituents are their twists. A variant that certifies only ellipticity, through exact rational linear-programming duals and without full rank, also succeeds on all $44$ blocks.

The two methods agree. In particular, the Frobenius traces of $J(X_\Gamma)$ obtained from the decomposition of Section~\ref{sec:reduction} coincide with those of an independent point count on the model $X_{H_{120}}$ of Remark~\ref{rem2.4} (Table~\ref{tab:traces}).

\begin{table}[ht]
\centering\small
\caption{Frobenius traces $\Tr(\mathrm{Frob}_p\mid H^1)$ of $J(X_\Gamma)$, obtained both from the decomposition of Section \ref{sec:reduction} and from a direct point count on $X_{H_{120}}$. The contribution of $J_{\chi_5}$ vanishes for $p\equiv\pm2\pmod5$.}\label{tab:traces}
\begin{tabular}{l*{11}{r}}
\toprule
$p$ & $7$ & $11$ & $13$ & $17$ & $19$ & $23$ & $29$ & $31$ & $37$ & $41$ & $43$\\
\midrule
$J_0(4800)$ & $-8$ & $-36$ & $-18$ & $-30$ & $4$ & $8$ & $-2$ & $-16$ & $6$ & $-102$ & $-36$\\
$J_{\chi_5}$ & $0$ & $-48$ & $0$ & $0$ & $-16$ & $0$ & $-32$ & $-48$ & $0$ & $-144$ & $0$\\
\midrule
$J(X_\Gamma)$ & $-8$ & $-84$ & $-18$ & $-30$ & $-12$ & $8$ & $-34$ & $-64$ & $6$ & $-246$ & $-36$\\
\bottomrule
\end{tabular}
\end{table}

\appendix
\section{Newspaces of level dividing 4800}\label{app:spaces}

Table~\ref{tab:spaces} lists, for each $M\mid4800$ with a nonzero newspace of character $1$ or $\chi_5$, the dimension of $S_2^\new(M,\varepsilon)$ and the number of Galois orbits, split by the degree of the coefficient field: ``$a{+}b$'' means $a$ rational and $b$ quadratic orbits for $\varepsilon=1$, and $a$ quadratic and $b$ quartic orbits for $\varepsilon=\chi_5$. An entry ``$0$, --'' means the space is zero; an empty entry means $5\nmid M$. The column $\sigma_0$ gives the multiplicity $\sigma_0(4800/M)$. The weighted sums are $\sum\sigma_0\dim=889$ for $\varepsilon=1$ and $888$ for $\varepsilon=\chi_5$. Divisors $M$ not listed have both spaces zero.

\begin{table}[h!]
\centering\small
\caption{Newspaces $S_2^\new(M,\varepsilon)$, $M\mid4800$, $\varepsilon\in\{1,\chi_5\}$.}\label{tab:spaces}
\begin{tabular}{rr rr rr}
\toprule
& & \multicolumn{2}{c}{$\varepsilon=1$} & \multicolumn{2}{c}{$\varepsilon=\chi_5$}\\
\cmidrule(lr){3-4}\cmidrule(lr){5-6}
$M$ & $\sigma_0$ & dim & orbits & dim & orbits\\
\midrule
15 & 14 & 1 & 1 & 0 & -- \\
20 & 20 & 1 & 1 & 0 & -- \\
24 & 12 & 1 & 1 & \multicolumn{2}{c}{} \\
30 & 12 & 1 & 1 & 2 & 1 \\
32 & 12 & 1 & 1 & \multicolumn{2}{c}{} \\
40 & 16 & 1 & 1 & 2 & 1 \\
48 & 9 & 1 & 1 & \multicolumn{2}{c}{} \\
50 & 12 & 2 & 2 & 2 & 1 \\
60 & 10 & 0 & -- & 2 & 1 \\
64 & 6 & 1 & 1 & \multicolumn{2}{c}{} \\
75 & 7 & 3 & 3 & 4 & 2 \\
80 & 12 & 2 & 2 & 2 & 1 \\
96 & 6 & 2 & 2 & \multicolumn{2}{c}{} \\
100 & 10 & 1 & 1 & 2 & 1 \\
120 & 8 & 2 & 2 & 2 & 1 \\
150 & 6 & 3 & 3 & 2 & 1 \\
160 & 8 & 4 & 2+1 & 6 & 1+1 \\
192 & 3 & 4 & 4 & \multicolumn{2}{c}{} \\
200 & 8 & 5 & 5 & 4 & 2 \\
240 & 6 & 4 & 4 & 6 & 3 \\
300 & 5 & 4 & 4 & 2 & 1 \\
320 & 4 & 8 & 6+1 & 10 & 3+1 \\
400 & 6 & 8 & 8 & 8 & 4 \\
480 & 4 & 8 & 8 & 12 & 4+1 \\
600 & 4 & 9 & 9 & 10 & 5 \\
800 & 4 & 19 & 9+5 & 18 & 5+2 \\
960 & 2 & 16 & 16 & 24 & 10+1 \\
1200 & 3 & 19 & 19 & 18 & 9 \\
1600 & 2 & 35 & 25+5 & 34 & 13+2 \\
2400 & 2 & 38 & 34+2 & 36 & 18 \\
4800 & 1 & 76 & 72+2 & 72 & 36 \\
\bottomrule
\end{tabular}
\end{table}

\section{Reproducibility}\label{app:repro}
\sloppy

The supplementary files that were used for the calculations in this article can be found at
\begin{center}
\url{https://github.com/rfauffar/genus1777}.
\end{center}

\subsection{The proof of Sections~\ref{sec:reduction}--\ref{sec:trivial}}
The script \texttt{verify\_genus1777.py} of the supplementary material uses PARI/GP \cite{PARI} (version 2.17.2, through the \texttt{cypari2} interface) and performs the following exact checks:
\begin{enumerate}[label=(V\arabic*),start=0]
\item all $70$ newspace dimensions $\dim S_2^\new(M,\varepsilon)$, $M\mid4800$, $\varepsilon\in\{1,\chi_5\}$ (with $5\mid M$ for $\chi_5$), agree with an independent evaluation, in plain Python without PARI, of the Cohen--Oesterl\'e dimension formula followed by inversion of $\dim S_2(N,\varepsilon)=\sum_M\sigma_0(N/M)\dim S_2^\new(M,\varepsilon)$;
\item $\dim S_2(4800,1)=889$ and $\dim S_2(4800,\chi_5)=888$, and the same total $1777$ is recovered as $\sum\sigma_0(4800/M)[K_f:\Q]$ over all newform orbits;
\item the decomposition of every space $S_2^\new(M,\varepsilon)$ in Table~\ref{tab:spaces} into Galois orbits, together with the assertion that the degrees of the orbits add up to the dimension of the space, so that the list of orbits is complete; in particular the list of the sixteen non-rational orbits of trivial character (computation (C1));
\item for each of the sixteen orbits, the coefficient identity of Lemma~\ref{lemma2} for every $n\le B(L)$ with $\gcd(n,D)=1$, in exact arithmetic in $K_f$ (computation (C2)); as a redundant check of the norm identities of Section~\ref{sec:trivial}, the script also evaluates the Hilbert symbols $(m,D)_v$ exactly at $v=\infty,2,3,5$, the only places where they can be nontrivial;
\item a consistency check of Proposition~\ref{prop:quartic}: for every Galois orbit of newforms of character $\chi_5$ and level dividing $4800$, the numbers $\psi^{-1}(p)a_p(f)$, $7\leq p<400$, are the eigenvalues of a newform of trivial character and level dividing $4800$;
\item the traces in Table~\ref{tab:traces}.
\end{enumerate}
The complete output of the script is included in the supplementary material.

A second script, \texttt{crosscheck\_modsym.py}, recomputes computation (C1) with PARI's modular-symbols package, which uses an algorithm independent of the one in (V2). For every $M\mid4800$ it splits the new subspace of the weight-two modular symbols for $\Gamma_0(M)$ (sign $+1$) into Hecke-irreducible pieces and checks that:
\begin{enumerate}[label=(X\arabic*)]
\item the dimensions of the pieces coincide with the degrees of the newform orbits from (V2);
\item for $p=7,11,13$, the characteristic polynomial of $T_p$ on the new subspace equals the product over the orbits of the characteristic polynomials of $a_p(f)$;
\item every piece of dimension $>1$ has the same characteristic polynomials of $T_7,T_{11},T_{13},T_{17}$ as exactly one of the sixteen orbits of Table~\ref{tab:sixteen}.
\end{enumerate}
For readers with access to Magma, the supplementary material also contains a script, \texttt{verify\_genus1777.m}, that repeats computations (C1) and (C2) in Magma. For levels dividing $2400$, the dimensions and orbit counts in Table~\ref{tab:spaces} also agree with those obtained independently in Magma \cite{Magma} for the verification of \S\ref{sec:equivariant}.

\subsection{The equivariant certificate}
The certificate of \S\ref{sec:equivariant} is contained in the compressed file
\begin{center}
\texttt{cert1777\_full.json.gz}
\end{center}
of the supplementary material. Its SHA-256 digest is
\begin{center}
\texttt{7cb552e4d3ba8e03afcef1323df74e82}\\ \texttt{4e2ff80673e6bb9317e9f4669c4c42c6}.
\end{center}
It contains the local character tables, the subgroup elements, the $44$ blocks, the candidate spaces, exact coefficient data, the $15576$ exact traces, the multiplicities and maximal minors, and the coefficient fields and twist data of the non-rational forms. Python verifiers recheck the group computations, the trace formula, the candidate lists and the rank minors, and the candidate eigenvalues were recomputed independently in Magma and compared with the data of \cite{PSdata} at all $354$ primes.


\begin{thebibliography}{99}

\bibitem{AL70}
A.~O.~L. Atkin and J.~Lehner, \emph{Hecke operators on $\Gamma_0(m)$}, Math. Ann. \textbf{185} (1970), 134--160.

\bibitem{AL78}
A.~O.~L. Atkin and W.-C.~W. Li, \emph{Twists of newforms and pseudo-eigenvalues of $W$-operators}, Invent. Math. \textbf{48} (1978), 221--243.

\bibitem{BBB+}
A.~J. Best, J.~Bober, A.~R. Booker, E.~Costa, J.~Cremona, M.~Derickx, M.~Lee, D.~Lowry-Duda, D.~Roe, A.~V. Sutherland, and J.~Voight, \emph{Computing classical modular forms}, in: Arithmetic Geometry, Number Theory, and Computation, Simons Symp., Springer, Cham, 2021; arXiv:2002.04717.

\bibitem{Magma}
W.~Bosma, J.~Cannon, and C.~Playoust, \emph{The Magma algebra system. I. The user language}, J. Symbolic Comput. \textbf{24} (1997), 235--265.

\bibitem{Cas73}
W.~Casselman, \emph{On some results of Atkin and Lehner}, Math. Ann. \textbf{201} (1973), 301--314.

\bibitem{ChenLuZuo}
K.~Chen, X.~Lu, and K.~Zuo, \emph{On CM points away from the Torelli locus}, J. London Math. Soc. \textbf{104} (2021), 1363--1383.

\bibitem{ChenLuZuo2}
K.~Chen, X.~Lu, and K.~Zuo, \emph{On a question of Ekedahl and Serre}, preprint, arXiv:1703.10377.

\bibitem{CW34}
C.~Chevalley and A.~Weil, \emph{\"Uber das Verhalten der Integrale 1. Gattung bei Automorphismen des Funktionenk\"orpers}, Abh. Math. Sem. Univ. Hamburg \textbf{10} (1934), 358--361.

\bibitem{DS05}
F.~Diamond and J.~Shurman, \emph{A First Course in Modular Forms}, Grad. Texts in Math. \textbf{228}, Springer, New York, 2005.

\bibitem{DT02}
W.~Duke and \'A.~T\'oth, \emph{The splitting of primes in division fields of elliptic curves}, Experiment. Math. \textbf{11} (2002), 555--565.

\bibitem{ES93}
T.~Ekedahl and J.-P.~Serre, \emph{Exemples de courbes alg\'ebriques \`a jacobienne compl\`etement d\'ecomposable}, C. R. Acad. Sci. Paris S\'er. I Math. \textbf{317} (1993), 509--513.

\bibitem{Fal83}
G.~Faltings, \emph{Endlichkeitss\"atze f\"ur abelsche Variet\"aten \"uber Zahlk\"orpern}, Invent. Math. \textbf{73} (1983), 349--366.


\bibitem{Li75}
W.-C.~W. Li, \emph{Newforms and functional equations}, Math. Ann. \textbf{212} (1975), 285--315.

\bibitem{Mom81}
F.~Momose, \emph{On the $\ell$-adic representations attached to modular forms}, J. Fac. Sci. Univ. Tokyo Sect. IA Math. \textbf{28} (1981), 89--109.

\bibitem{PARI}
The PARI~Group, \emph{PARI/GP version 2.17.2}, Univ. Bordeaux, \url{https://pari.math.u-bordeaux.fr/}.

\bibitem{PR17}
J.~Paulhus and A.~M. Rojas, \emph{Completely decomposable Jacobian varieties in new genera}, Exp. Math. \textbf{26} (2017), no.~4, 430--445.

\bibitem{PS25}
J.~Paulhus and A.~V. Sutherland, \emph{Completely decomposable modular Jacobians}, Contemp. Math. \textbf{840} (2026), 261--270.

\bibitem{PSdata}
J.~Paulhus and A.~V. Sutherland, \emph{Code and data for ``Completely decomposable modular Jacobians''}, \url{https://github.com/AndrewVSutherland/CompletelyDecomposableModularJacobians}.

\bibitem{Pyl04}
E.~E. Pyle, \emph{Abelian varieties over $\Q$ with large endomorphism algebras and their simple components over $\Qbar$}, in: Modular Curves and Abelian Varieties, Progr. Math. \textbf{224}, Birkh\"auser, Basel, 2004, 189--239.

\bibitem{Que09}
J.~Quer, \emph{Fields of definition of building blocks}, Math. Comp. \textbf{78} (2009), 537--554.

\bibitem{Rib80}
K.~A. Ribet, \emph{Twists of modular forms and endomorphisms of abelian varieties}, Math. Ann. \textbf{253} (1980), 43--62.

\bibitem{Rib92}
K.~A. Ribet, \emph{Abelian varieties over $\Q$ and modular forms}, in: Algebra and Topology 1992 (Taej\u{o}n), Korea Adv. Inst. Sci. Tech., Taej\u{o}n, 1992, 53--79; reprinted in Progr. Math. \textbf{224}, Birkh\"auser, 2004, 241--261.

\bibitem{RodRoj}
R.~E. Rodr\'{\i}guez and A.~M. Rojas, \emph{Decomposing abelian varieties into simple factors: algorithms and applications}, J. Algebra \textbf{683} (2025), 364--388.

\bibitem{RSZB}
J.~Rouse, A.~V. Sutherland, and D.~Zureick-Brown, \emph{$\ell$-adic images of Galois for elliptic curves over $\Q$}, with an appendix by J.~Voight, Forum Math. Sigma \textbf{10} (2022), e62; arXiv:2106.11141.

\bibitem{Ser20}
J.-P. Serre, \emph{Rational points on curves over finite fields}, with contributions by E.~Howe, J.~Oesterl\'e and C.~Ritzenthaler, Documents Math\'ematiques \textbf{18}, Soc. Math. France, Paris, 2020.

\bibitem{Shi71}
G.~Shimura, \emph{Introduction to the Arithmetic Theory of Automorphic Functions}, Publ. Math. Soc. Japan \textbf{11}, Iwanami Shoten and Princeton Univ. Press, 1971.

\bibitem{Shi71b}
G.~Shimura, \emph{On elliptic curves with complex multiplication as factors of the Jacobians of modular function fields}, Nagoya Math. J. \textbf{43} (1971), 199--208.

\bibitem{Stu87}
J.~Sturm, \emph{On the congruence of modular forms}, in: Number Theory (New York, 1984--1985), Lecture Notes in Math. \textbf{1240}, Springer, Berlin, 1987, 275--280.

\bibitem{Yam}
T.~Yamauchi, \emph{On $\Q$-simple factors of Jacobian varieties of modular curves}, Yokohama Math. J. \textbf{53} (2007), no.~2, 149--160.

\end{thebibliography}
\end{document}